\documentclass[10pt]{amsart}
\pdfoutput=1
\usepackage{amscd}
\usepackage{amsmath}
\usepackage{amssymb}
\usepackage{graphicx}
\usepackage[all, pdf]{xy}
\usepackage{hyperref}

\usepackage[T1]{fontenc}

\title[Derived Categories]
{Introduction to Derived Categories}

\author{Amnon Yekutieli}

\address{Department of  Mathematics,
Ben Gurion University, Be'er Sheva 84105, Israel}
\email{amyekut@math.bgu.ac.il}

\date{27 Jan 2015}

\newtheorem{thm}[equation]{Theorem}
\newtheorem{cor}[equation]{Corollary}
\newtheorem{prop}[equation]{Proposition}

\theoremstyle{definition}
\newtheorem{dfn}[equation]{Definition}

\newtheorem{exa}[equation]{Example}

\newtheorem{exer}[equation]{Exercise}

\numberwithin{equation}{section}

\newcommand{\iso}{\xrightarrow{\simeq}}

\newcommand{\xar}{\xrightarrow}
\newcommand{\opn}{\operatorname}
\newcommand{\cat}[1]{\operatorname{\mathsf{#1}}}

\newcommand{\rmitem}[1]{\item[\text{\textup{(#1)}}]}
\newcommand{\mfrak}[1]{\mathfrak{#1}}
\newcommand{\mcal}[1]{\mathcal{#1}}

\newcommand{\mrm}[1]{\mathrm{#1}}
\newcommand{\mbb}[1]{\mathbb{#1}}

\newcommand{\tup}[1]{\textup{#1}}
\newcommand{\bsym}[1]{\boldsymbol{#1}}

\newcommand{\ot}{\otimes}
\newcommand{\til}[1]{\tilde{#1}}

\newcommand{\K}{\mbb{K}}

\newcommand{\N}{\mbb{N}}
\newcommand{\Z}{\mbb{Z}}

\newcommand{\g}{\mfrak{g}}
\newcommand{\m}{\mfrak{m}}

\newcommand{\p}{\mfrak{p}}

\newcommand{\bwedge}{{\textstyle \bigwedge}}

\renewcommand{\d}{\mrm{d}}

\newcommand{\ga}{\gamma}
\newcommand{\al}{\alpha}
\newcommand{\be}{\beta}

\newcommand{\OO}{\mcal{O}}
\newcommand{\MM}{\mcal{M}}

\newcommand{\II}{\mcal{I}}

\renewcommand{\AA}{\mcal{A}}

\renewcommand{\th}{\theta}

\newcommand{\dcat}[1]{\operatorname{\bsym{\cat{#1}}}}

\newcommand{\lb}{\linebreak}
\newcommand{\bmat}[1]{\begin{bmatrix} #1 \end{bmatrix}}
\newcommand{\cd}{\,{\cdot}\,}

\begin{document}

\begin{abstract}
Derived categories were invented by Grothendieck and Verdier around 1960,
not very long after the ``old'' homological algebra (of derived functors 
between abelian categories) was established. This ``new'' homological algebra, 
of derived categories and derived functors between them, provides a 
significantly richer and more flexible machinery than the ``old'' homological
algebra. For instance, the important concepts of {\em dualizing complex} and 
{\em tilting complex} do not exist in the ``old'' homological algebra. 

This paper is an edited version of the notes for a two-lecture minicourse 
given at MSRI in January 2013. Sections 1-5 are about the general theory of 
derived categories, and the material is taken from my manuscript ``A Course on 
Derived Categories'' (available online). Sections 6-9 are on more specialized 
topics, leaning towards noncommutative algebraic geometry.  
\end{abstract}
 
\maketitle

\tableofcontents

\setlength{\parskip}{1ex} 
\setlength{\parindent}{0ex}

\section{The Homotopy Category}

Suppose $\cat{M}$ is an {\em  abelian category}. 
The main examples for us are these:

\begin{itemize}
\item $A$ is a ring, and 
$\cat{M} = \cat{Mod} A$, the category of left $A$-modules. 

\item $(X, \AA)$ is a ringed space, and 
$\cat{M} = \cat{Mod} \AA$, the category of sheaves of left $\AA$-modules.
\end{itemize}

A {\em complex} in $\cat{M}$ is a diagram 
\[ M = \bigl( \cdots \to M^{-1} \xar{\ \d_M^{-1} \ }  M^{0} \xar{\ \d_M^0 \ } 
 M^{1} \to \cdots \bigr)  \]
in $\cat{M}$ such that $\d_M^{i+1} \circ \d_M^i = 0$. 
A morphism of complexes 
$\phi : M \to N$  is a commutative diagram 
\begin{equation}  \label{eqn:1}
\UseTips \xymatrix @C=5ex @R=6ex {
\cdots 
\ar[r]
& 
M^{-1}
\ar[r]^{\d_M^{-1}}
\ar[d]^{\phi^{-1}}
&
M^{0} 
\ar[r]^{\d_M^{0}}
\ar[d]^{\phi^{0}}
&
M^{1}
\ar[d]^{\phi^{1}}
\ar[r]
& 
\cdots
\\
\cdots 
\ar[r]
& 
N^{-1} 
\ar[r]_{\d_N^{-1}}
&
N^{0} 
\ar[r]_{\d_N^{0}}
&
N^{1} 
\ar[r]
&
\cdots
}
\end{equation}
in  $\cat{M}$.
Let us denote by {\em $\dcat{C}(\cat{M})$ the category of complexes} in
$\cat{M}$.
It is again an abelian category; but it is also a {\em differential graded 
category}, as we now explain.  

Given $M, N \in \dcat{C}(\cat{M})$ we let 
\[ \opn{Hom}_{\cat{M}}(M, N)^i := 
\prod_{j \in \Z} \opn{Hom}_{\cat{M}}(M^j, N^{j + i}) \]
and 
\[ \opn{Hom}_{\cat{M}}(M, N) := 
\bigoplus_{i \in \Z} \opn{Hom}_{\cat{M}}(M, N)^i . \]
For $\phi \in \opn{Hom}_{\cat{M}}(M, N)^i$ we let 
\[ \d(\phi) := \d_N \circ \phi - (-1)^i \cd \phi \circ \d_M . \]
In this way $\opn{Hom}_{\cat{M}}(M, N)$ becomes a complex of abelian
groups, i.e.\ a DG (differential graded) $\Z$-module. 
Given a third complex $L \in \dcat{C}(\cat{M})$, composition of morphisms in 
$\cat{M}$ induces a homomorphism of DG $\Z$-modules 
\[ \opn{Hom}_{\cat{M}}(L, M) \ot_{\Z} \opn{Hom}_{\cat{M}}(M, N) \to 
\opn{Hom}_{\cat{M}}(L, N) . \]
Cf.\ Section 5; a DG algebra is a DG category with one object. 

Note that the abelian structure of $\dcat{C}(\cat{M})$ can be recovered from the
DG structure  as follows:
\[ \opn{Hom}_{\dcat{C}(\cat{M})}(M, N) = 
\opn{Z}^0 \bigl( \opn{Hom}_{\cat{M}}(M, N) \bigr) , \]
the set of $0$-cocycles. 
Indeed, for $\phi : M \to N$ of degree $0$
the condition $\d(\phi) = 0$ 
is equivalent to the commutativity of the diagram (\ref{eqn:1}).

Next we define the {\em homotopy category $\dcat{K}(\cat{M})$}.
Its objects are the complexes in $\cat{M}$ (same as $\dcat{C}(\cat{M})$), 
and 
\[ \opn{Hom}_{\dcat{K}(\cat{M})}(M, N) = 
\opn{H}^0 \bigl( \opn{Hom}_{\cat{M}}(M, N) \bigr) . \]
In other words, these are homotopy classes of morphisms 
$\phi : M \to N$ in $\dcat{C}(\cat{M})$.
 
There is an additive functor 
$\dcat{C}(\cat{M}) \to \dcat{K}(\cat{M})$,
which is the identity on objects and surjective on morphisms. 
The additive category $\dcat{K}(\cat{M})$ is no longer abelian -- it is a
{\em triangulated category}. Let me explain what this means.

Suppose $\cat{K}$ is an additive category, with an automorphism $\opn{T}$ 
called the {\em translation} (or  shift, or  suspension).
A {\em triangle} in $\cat{K}$ is a diagram of morphisms of this sort:
\[ L \xar{\al} M \xar{\be} N \xar{\ga} \opn{T}(L) . \]
The name comes from the alternative typesetting
\[ \UseTips  \xymatrix @C=4ex @R=4ex {
& 
N 
\ar[dl]_{\ga}
\\
L
\ar[rr]^{\al}
& &
M
\ar[lu]_{\be}
} \]
A triangulated category structure on $\cat{K}$ is a set of triangles called
{\em distinguished triangles}, satisfying a list of axioms (that are not so
important for us). 
Details can be found in the references  \cite{Ye6}, \cite{Sc}, 
\cite{Ha}, \cite{We}, \cite{KS1}, \cite{Ne2} or \cite{LH}. 

The translation $\opn{T}$ of the category $\dcat{K}(\cat{M})$ is defined as 
follows. On objects we take
$\opn{T}(M)^i := M^{i+1}$ and $\d_{\opn{T}(M)} := - \d_M$. On morphisms it is 
$\opn{T}(\phi)^i := \phi^{i+1}$. 
For $k \in \Z$, the $k$-th translation of $M$ is denoted by 
$M[k] := \opn{T}^k(M)$.
 
Given a morphism $\al : L \to M$ in $\dcat{C}(\cat{M})$, its {\em cone}
is the complex 
\[  \opn{cone}(\al) := \opn{T}(L) \oplus M = \bmat{\opn{T}(L) \\ M} \]
with differential (in matrix notation)
\[ \d := 
\bmat{ \d_{\opn{T}(L)} & 0 \\[0.5em] \al &  \d_M } ,  \]
where $\al$ is viewed as a degree $1$ morphism $\opn{T}(L) \to M$. 
There are canonical morphisms 
$M \to \opn{cone}(\al)$ and $\opn{cone}(\al) \to \opn{T}(L)$ in 
$\dcat{C}(\cat{M})$.

A triangle in $\dcat{K}(\cat{M})$ is distinguished 
if it is  isomorphic, as a diagram in $\dcat{K}(\cat{M})$, to the triangle 
\[ L \xar{\al} M \xar{} \opn{cone}(\al) \xar{} \opn{T}(L) \]
for some morphism $\al : L \to M$ in $\dcat{C}(\cat{M})$.
A calculation shows that $\dcat{K}(\cat{M})$ is indeed triangulated (i.e.\ the
axioms that I did not specify are satisfied).
 
The relation between distinguished triangles and exact
sequences  will be mentioned later.

Suppose $\cat{K}$ and $\cat{K}'$ are triangulated categories. A
{\em triangulated functor} $F : \cat{K} \to \cat{K}'$ is an additive functor
that commutes with the translations, and sends distinguished triangles to
distinguished triangles.
 
\begin{exa} \label{exa:1}
Let $F : \cat{M} \to \cat{M}'$ be an additive functor
(not necessarily exact) between abelian categories. 
Extend $F$ to a functor 
\[ \dcat{C}(F) : \dcat{C}(\cat{M}) \to  \dcat{C}(\cat{M}') \]
in the obvious way, namely 
\[ \dcat{C}(F)(M)^i := F(M^i) \]
for a complex $M = \{ M^i \}_{i \in \Z}$.
The functor $\dcat{C}(F)$ respects homotopies, so we get an additive functor 
\[ \dcat{K}(F) : \dcat{K}(\cat{M}) \to  \dcat{K}(\cat{M}') . \]
This is a triangulated functor. 
\end{exa}

\section{The Derived Category}

As before $\cat{M}$ is an abelian category.
Given a complex $M \in \dcat{C}(\cat{M})$, we can consider its cohomologies 
\[ \mrm{H}^i(M) := \opn{ker}(\d_M^i) /  \opn{im}(\d_M^{i-1}) \in  \cat{M} . \]
Since the cohomologies are homotopy-invariant, we get additive functors 
\[ \mrm{H}^i : \dcat{K}(\cat{M}) \to  \cat{M} . \]
 
A morphism $\psi : M \to N$ in $\dcat{K}(\cat{M})$ is called a 
{\em quasi-isomorphism} if $\mrm{H}^i(\psi)$ are isomorphisms for all $i$. 
Let us denote by $\dcat{S}(\cat{M})$ the set of all quasi-isomorphisms in 
$\dcat{K}(\cat{M})$. 
Clearly $\dcat{S}(\cat{M})$ is a multiplicatively closed set, i.e.\ the
composition of two quasi-isomorphisms is a quasi-isomorphism. 
A calculation shows that $\dcat{S}(\cat{M})$ is a left and right
denominator set (as in ring theory). 
It follows that the Ore localization $\dcat{K}(\cat{M})_{\dcat{S}(\cat{M})}$
exists. This is an additive category, with object set 
\[ \opn{Ob}(\dcat{K}(\cat{M})_{\dcat{S}(\cat{M})}) = 
\opn{Ob}(\dcat{K}(\cat{M})) . \]
There is a functor 
\[ \opn{Q} : \dcat{K}(\cat{M}) \to \dcat{K}(\cat{M})_{\dcat{S}(\cat{M})} \]
called the localization functor, which is the identity on objects.
Every morphism 
$\chi : M \to N$ in $\dcat{K}(\cat{M})_{\dcat{S}(\cat{M})}$ can be written as 
\[ \chi = \opn{Q}(\phi_1) \circ \opn{Q}(\psi_1^{-1}) = 
\opn{Q}(\psi_2^{-1}) \circ \opn{Q}(\phi_2) \]
for some $\phi_i \in \dcat{K}(\cat{M})$
and $\psi_i \in \dcat{S}(\cat{M})$.
 
The category $\dcat{K}(\cat{M})_{\dcat{S}(\cat{M})}$ inherits a triangulated
structure from $\dcat{K}(\cat{M})$, 
and the localization functor $\opn{Q}$ is triangulated.
There is a universal property: given a triangulated functor 
\[ F : \dcat{K}(\cat{M}) \to \cat{E} \]
to a triangulated category $\cat{E}$, such that $F(\psi)$ is an isomorphism for
every $\psi \in \dcat{S}(\cat{M})$, there exists a unique triangulated functor 
\[ F_{\dcat{S}(\cat{M})} : \dcat{K}(\cat{M})_{\dcat{S}(\cat{M})} \to \cat{E}  \]
such that 
\[ F_{\dcat{S}(\cat{M})}  \circ \opn{Q} = F . \]

\begin{dfn}
The {\em derived category} of the abelian category $\cat{M}$ is the
triangulated
category 
\[ \dcat{D}(\cat{M}) := \dcat{K}(\cat{M})_{\dcat{S}(\cat{M})} . \]
\end{dfn}
 
The derived category was introduced by Grothendieck and Verdier around 1960. 
The first published material is the book ``Residues and Duality'' \cite{Ha}
from 1966, written by Hartshorne following notes by Grothendieck. 
 
Let $\dcat{D}(\cat{M})^0$ be the full subcategory of $\dcat{D}(\cat{M})$
consisting of the complexes whose cohomology is concentrated in degree $0$.
 
\begin{prop}
The obvious functor $\cat{M} \to \dcat{D}(\cat{M})^0$ is an equivalence.
\end{prop}
 
This allows us to view $\cat{M}$ as an additive subcategory of 
$\dcat{D}(\cat{M})$.
It turns out that the abelian structure of $\cat{M}$ can be
recovered from this embedding. 
 
\begin{prop}
Consider a sequence
\[ 0 \to L \xar{\al} M \xar{\be} N \to 0 \]
in $\cat{M}$. 
This sequence is exact iff there is a morphism 
$\ga : N \to L[1]$ in $\dcat{D}(\cat{M})$ such that 
\[  L \xar{\al} M \xar{\be} N \xar{\ga} L[1] \]
is a distinguished triangle.
\end{prop}

\section{Derived Functors}

As before $\cat{M}$ is an abelian category. Recall the localization functor 
\[ \opn{Q} : \dcat{K}(\cat{M}) \to \dcat{D}(\cat{M}) . \]
It is a triangulated functor, 
which is the identity on objects, and inverts quasi-isomorphisms.
 
Suppose $\cat{E}$ is some triangulated category, and 
$F : \dcat{K}(\cat{M}) \to \cat{E}$ a triangulated functor.
We now introduce the right and left derived functors of $F$. These are
triangulated functors
\[ \mrm{R} F, \mrm{L} F : \dcat{D}(\cat{M}) \to \cat{E}  \]
satisfying suitable universal properties. 

\begin{dfn} 
A {\em right derived functor} of $F$ 
is a triangulated functor 
\[ \mrm{R} F : \dcat{D}(\cat{M}) \to \cat{E} , \]
together with a morphism 
\[ \eta : F \to \mrm{R} F \circ \opn{Q} \]
of triangulated functors $\dcat{K}(\cat{M}) \to \cat{E}$,
satisfying this condition: 
\begin{itemize} 
\item[($*$)] The pair $(\mrm{R} F, \eta)$ is initial among all such pairs.
\end{itemize}
\end{dfn}

Being initial means that if $(G, \eta')$ is another such pair,
then  there is a unique morphism of triangulated functors 
$\th :  \mrm{R} F \to G$
s.t.\ $\eta' = \th \circ \eta$. 
The universal condition implies that if a right derived functor 
$(\mrm{R} F, \eta)$ exists, then it is unique, up to
a unique isomorphism of triangulated functors.
 
\begin{dfn}
A {\em left derived functor} of $F$ 
is a triangulated functor 
\[ \mrm{L} F : \dcat{D}(\cat{M}) \to \cat{E} , \]
together with a morphism 
\[ \eta : \mrm{L} F \circ \opn{Q} \to F  \]
of triangulated functors $\dcat{K}(\cat{M}) \to \cat{E}$,
satisfying this condition: 
\begin{itemize} 
\item[($*$)] The pair $(\mrm{L} F, \eta)$ is terminal among all such pairs.
\end{itemize}
\end{dfn}

Again, if $(\mrm{L} F, \eta)$ exists, then it is unique up to a unique
isomorphism.

There are various modifications. 
One of them is a contravariant triangulated functor
$F  : \dcat{K}(\cat{M}) \to \cat{E}$. 
This can be handled using the fact that
$\dcat{K}(\cat{M})^{\mrm{op}}$ is triangulated, and 
$F  : \dcat{K}(\cat{M})^{\mrm{op}} \to \cat{E}$ is covariant. 
 
We will also want to derive bifunctors. 
Namely to a bitriangulated bifunctor 
\[ F  : \dcat{K}(\cat{M}) \times \dcat{K}(\cat{M}')  \to \cat{E} \]
we will want to associate bitriangulated bifunctors
\[ \mrm{R} F,  \mrm{L} F  : \dcat{D}(\cat{M}) \times \dcat{D}(\cat{M}')  \to
\cat{E} . \]
This is done similarly, and I won't give details.

\section{Resolutions}
  
Consider an additive functor 
$F : \cat{M} \to \cat{M}'$ between abelian categories, and the corresponding 
triangulated functor 
$\dcat{K}(F) : \dcat{K}(\cat{M}) \to  \dcat{K}(\cat{M}')$, as in Example
\ref{exa:1}.
By slight abuse we write $F$ instead of $\dcat{K}(F)$. 
We want to construct (or prove existence) of the derived functors 
\[ \mrm{R} F, \mrm{L} F : \dcat{D}(\cat{M}) \to \dcat{D}(\cat{M}') . \]
If $F$ is exact (namely $F$ sends exact sequences to exact sequences), then 
$\mrm{R} F = \mrm{L} F = F$. (This is an easy
exercise.) Otherwise we need {\em resolutions}. 

The DG structure of $\dcat{C}(\cat{M})$ gives, for every 
$M, N \in \dcat{C}(\cat{M})$, a complex of abelian groups
$\opn{Hom}_{\cat{M}}(M, N)$. 
Recall that a complex $N$ is called acyclic if $\mrm{H}^i(N) = 0$ for all $i$; 
i.e.\ $N$ is an exact sequence in $\cat{M}$. 
 
\begin{dfn}
\begin{enumerate} 
\item A complex $I \in \dcat{K}(\cat{M})$ is called {\em K-injective} if for
every acyclic $N \in \dcat{K}(\cat{M})$, the complex 
$\opn{Hom}_{\cat{M}}(N, I)$ is also acyclic.
 
\item Let $M \in \dcat{K}(\cat{M})$. A {\em K-injective resolution} of $M$ is
a quasi-isomorphism $M \to I$ in $\dcat{K}(\cat{M})$, where $I$ is K-injective.
 
\item We say that $\dcat{K}(\cat{M})$ {\em has enough K-injectives} if every 
$M \in \dcat{K}(\cat{M})$ has some K-injective resolution.
\end{enumerate}
\end{dfn}

\begin{thm}
If $\dcat{K}(\cat{M})$ has enough K-injectives, then every 
triangulated functor 
$F : \dcat{K}(\cat{M}) \to \cat{E}$ 
has a right derived functor 
$(\mrm{R} F, \eta)$. 
Moreover, for every K-injective complex $I \in \dcat{K}(\cat{M})$, the
morphism $\eta_I : F(I) \to \mrm{R} F (I)$ 
in $\cat{E}$ is an isomorphism.
\end{thm}
 
The proof / construction goes like this: for every 
$M \in \dcat{K}(\cat{M})$ we choose a K-injective resolution 
$\zeta_M : M \to I_M$, and we define 
\[ \mrm{R} F (M) := F (I_M) \]
and
\[ \eta_M := F(\zeta_M) : F(M) \to F (I_M) \]
in $\cat{E}$.

Regarding existence of K-injective resolutions: 
 
\begin{prop}
A bounded below complex of injective objects of $\cat{M}$ is a K-injective
complex.
\end{prop}
 
This is the type of injective resolution used in \cite{Ha}. 
The most general statement I know is this 
(see \cite[Theorem 14.3.1]{KS2}): 
 
\begin{thm}
If $\cat{M}$ is a Grothendieck abelian category, then 
$\dcat{K}(\cat{M})$ has enough K-injectives.
\end{thm}
 
This includes  $\cat{M} = \cat{Mod} A$ for a ring $A$, 
and $\cat{M} = \cat{Mod} \AA$ for a sheaf of rings $\AA$.
Actually in these cases the construction of K-injective
resolutions can be done very explicitly, and it is not so difficult. 

\begin{exa}
Let $f : (X, \AA_X) \to (Y, \AA_Y)$ be a map of ringed spaces. 
(For instance a map of schemes 
$f : (X, \OO_X) \to (Y, \OO_Y)$.)
The map $f$ induces an additive functor 
\[ f_* : \cat{Mod} \AA_X \to \cat{Mod} \AA_Y  \]
called push-forward, which is usually not exact (it is left exact though). 
Since $\dcat{K}(\cat{Mod} \AA_X)$ has enough K-injectives, the right derived 
functor 
\[  \mrm{R} f_* : \dcat{D}(\cat{Mod} \AA_X) \to \dcat{D}(\cat{Mod} \AA_Y) \]
exists. 
 
For $\MM \in \cat{Mod} \AA_X$ we can use an injective
resolution $\MM \to \II$ (in the ``classical'' sense), and therefore 
\[ \mrm{H}^q(\mrm{R} f_*(\MM)) \cong \mrm{H}^q( f_*(\II)) \cong
\mrm{R}^q f_* (\MM) , \]
where the latter is the ``classical'' right derived functor. 
\end{exa}

Analogously we have:

\begin{dfn} 
\begin{enumerate} 
\item A complex $P \in \dcat{K}(\cat{M})$ is called {\em K-projective} if for
every acyclic $N \in \dcat{K}(\cat{M})$, the complex 
$\opn{Hom}_{\cat{M}}(P, N)$ is also acyclic.
 
\item Let $M \in \dcat{K}(\cat{M})$. A {\em K-projective resolution} of $M$
is a quasi-isomorphism $P \to M$ in $\dcat{K}(\cat{M})$, where $P$ is 
K-projective.
 
\item We say that $\dcat{K}(\cat{M})$ {\em has enough K-projectives} if
every  $M \in \dcat{K}(\cat{M})$ has some K-projective resolution.
\end{enumerate}
\end{dfn}
 
\begin{thm}
If $\dcat{K}(\cat{M})$ has enough K-projectives, then every 
triangulated functor 
$F : \dcat{K}(\cat{M}) \to \cat{E}$ 
has a left derived functor 
$(\mrm{L} F, \eta)$. 
Moreover, for every K-projective complex $P \in \dcat{K}(\cat{M})$, the
morphism $\eta_P : \mrm{L} F(P) \to F (P)$ 
in $\cat{E}$ is an isomorphism.
\end{thm}
 
The construction of $\mrm{L} F$ is by K-projective resolutions.

\begin{prop}
A bounded above complex of projective objects of $\cat{M}$ is a K-projective
complex.
\end{prop}
 
\begin{prop}
Let $A$ be a ring. Then $\dcat{K}(\cat{Mod} A)$ has enough K-projectives.
\end{prop}

The construction of K-projective resolutions in this case can be done very 
explicitly, and it is not difficult. 
 
The concepts of K-injective and K-projective complexes were introduced by
Spaltenstein \cite{Sp} in 1988. At about the same time other authors (Keller
\cite{Ke}, Bockstedt-Neeman \cite{BN}, \ldots) discovered these concepts
independently, with other names (such as {\em homotopically injective
complex}).

\begin{exa} \label{exa:5}
Suppose $\K$ is a commutative ring and $A$ is a $\K$-algebra (i.e.\ $A$ is a
ring and there is a homomorphism $\K \to \mrm{Z}(A)$). 
Consider the bi-additive bifunctor 
\[ \opn{Hom}_A(- , -) : (\cat{Mod} A)^{\mrm{op}} \times \cat{Mod} A 
\to \cat{Mod} \K . \]
 
We have seen how to extend this functor to complexes (this is sometimes called
``product totalization''), giving rise to a bitriangulated bifunctor 
\[ \opn{Hom}_A(- , -) : \dcat{K}(\cat{Mod} A)^{\mrm{op}} \times 
\dcat{K}(\cat{Mod} A) \to \dcat{K}(\cat{Mod} \K) . \]
The right derived bifunctor 
\[ \opn{RHom}_A(- , -) : \dcat{D}(\cat{Mod} A)^{\mrm{op}} \times 
\dcat{D}(\cat{Mod} A) \to \dcat{D}(\cat{Mod} \K)  \]
can be constructed / calculated by a K-injective resolution in either the
first or the second argument.
Namely given $M, N \in \dcat{K}(\cat{Mod} A)$ we can choose a K-injective 
resolution $N \to I$, and let 
\begin{equation}  \label{eqn:5}
\opn{RHom}_A(M , N) := \opn{Hom}_A(M , I) \in 
\dcat{D}(\cat{Mod} \K) .
\end{equation}
Or we can choose a K-injective resolution $M \to P$ in 
$\dcat{K}(\cat{Mod} A)^{\mrm{op}}$, which is really a K-projective resolution 
$P \to M$ in $\dcat{K}(\cat{Mod} A)$, and let
\begin{equation}  \label{eqn:6}
\opn{RHom}_A(M , N) := \opn{Hom}_A(P , N)
\in \dcat{D}(\cat{Mod} \K) .
\end{equation}
The two complexes (\ref{eqn:5}) and (\ref{eqn:6})  are canonically related by
the quasi-iso\-morphisms 
\[ \opn{Hom}_A(P , N) \to \opn{Hom}_A(P , I) \leftarrow \opn{Hom}_A(M , I) . \]
 
If $M, N \in \cat{Mod} A$ then of course
\[ \mrm{H}^q \bigl( \opn{RHom}_A(M , N) \bigr) \cong
\opn{Ext}^q_A(M , N) , \]
where the latter is ``classical'' Ext. 
\end{exa}

K-projective and K-injective complexes are good also for understanding the
structure of $\dcat{D}(\cat{M})$.
 
\begin{prop}
Suppose $P \in \dcat{K}(\cat{M})$ is K-projective and 
$I \in \dcat{K}(\cat{M})$ is K-injective. 
Then for any $M \in \dcat{K}(\cat{M})$
the homomorphisms 
\[ \opn{Q} : \opn{Hom}_{\dcat{K}(\cat{M})}(P , M) \to 
\opn{Hom}_{\dcat{D}(\cat{M})}(P , M) \]
and
\[ \opn{Q} : \opn{Hom}_{\dcat{K}(\cat{M})}(M , I) \to 
\opn{Hom}_{\dcat{D}(\cat{M})}(M , I) \]
are bijective. 
\end{prop}

Let us denote by $\dcat{K}(\cat{M})_{\tup{prj}}$ and 
$\dcat{K}(\cat{M})_{\tup{inj}}$ the full subcategories of 
$\dcat{K}(\cat{M})$ on the K-projective and the K-injective complexes
respectively. 
 
\begin{cor}
The triangulated functors 
\[ \opn{Q} : \dcat{K}(\cat{M})_{\tup{prj}} \to \dcat{D}(\cat{M}) \]
and 
\[ \opn{Q} : \dcat{K}(\cat{M})_{\tup{inj}} \to \dcat{D}(\cat{M}) \]
are fully faithful.
\end{cor}

\begin{cor}
\begin{enumerate}
\item If $\dcat{K}(\cat{M})$ has enough K-projectives, then the triangulated 
functor 
$\opn{Q} : \dcat{K}(\cat{M})_{\tup{prj}} \to \dcat{D}(\cat{M})$
is an equivalence. 

\item If $\dcat{K}(\cat{M})$ has enough K-injectives, then the triangulated 
functor
$\opn{Q} : \lb \dcat{K}(\cat{M})_{\tup{inj}} \to \dcat{D}(\cat{M})$
is an equivalence. 
\end{enumerate}
\end{cor}
 
\begin{exer}
Let $\K$ be a nonzero commutative ring and $A := \K[t]$ the polynomial
ring. We view $\K$ as an $A$-module via $t \mapsto 0$. Find a nonzero 
morphism  $\chi : \K \to \K[1]$ in $\dcat{D}(\cat{Mod} A)$.
Show that $\mrm{H}^q(\chi) = 0$ for all $q \in \Z$. 
\end{exer}

When $\cat{M} = \cat{Mod} A$ for a ring $A$, we can also talk about 
{\em K-flat} complexes. A complex $P$ is K-flat if for any acyclic complex 
$N \in \cat{Mod} A^{\mrm{op}}$ the complex
$N \ot_A P$ is acyclic. Any K-projective complex is K-flat. The left derived 
bifunctor $N \ot^{\mrm{L}}_A M$ can be constructed using K-flat resolutions of 
either argument: 
\[ N \ot^{\mrm{L}}_A M \cong N \ot_A P \cong Q \ot_A M \]
for any K-flat resolutions $P \to M$ in $\dcat{K}(\cat{Mod} A)$
and $Q \to N$ in $\dcat{K}(\cat{Mod} A^{\mrm{op}})$.

\section{DG Algebras}

A DG algebra (or DG ring) is a graded ring $A = \bigoplus_{i \in \Z} A^i$,
with differential $\d$ of degree $1$, satisfying the graded Leibniz rule 
\[ \d(a \cd b) = \d(a) \cd b + (-1)^i \cd a \cd \d(b) \]
for $a \in A^i$ and $b \in A^j$. 
 
A left DG $A$-module is a left graded $A$-module 
$M = \bigoplus_{i \in \Z} M^i$,
with differential $\d$ of degree $1$, satisfying the graded Leibniz rule.
Denote by $\cat{DGMod} A$ the category of left DG $A$-modules. 
 
As in the ring case, for any $M, N \in \cat{DGMod} A$  there is a complex of
$\Z$-modules 
$\opn{Hom}_A(M, N)$, and 
\[ \opn{Hom}_{\cat{DGMod} A}(M, N) = \mrm{Z}^0 \bigl( \opn{Hom}_A(M, N) \bigr).
\]
The homotopy category is $\til{\dcat{K}}(\cat{DGMod} A)$, with 
\[ \opn{Hom}_{\til{\dcat{K}}(\cat{DGMod} A)}(M, N) = 
\mrm{H}^0 \bigl( \opn{Hom}_A(M, N) \bigr) . \]
After inverting the quasi-isomorphisms in 
$\til{\dcat{K}}(\cat{DGMod} A)$ we obtain the derived category 
$\til{\dcat{D}}(\cat{DGMod} A)$. These are triangulated categories.
 
\begin{exa}
Suppose $A$ is a ring (i.e.\ $A^i = 0$ for $i \neq 0$). 
Then $\cat{DGMod} A = \dcat{C}(\cat{Mod} A)$ and 
$\til{\dcat{D}}(\cat{DGMod} A) = \dcat{D}(\cat{Mod} A)$.
\end{exa}

Derived functors are defined as in the ring case, and there are enough
K-injectives, K-projectives and K-flats in 
$\til{\dcat{K}}(\cat{DGMod} A)$. 
 
Let $A \to B$ be a homomorphism of DG algebras. There are 
additive functors 
\[ B \ot_A - : \
\cat{DGMod} A \rightleftarrows \cat{DGMod} B \
: \opn{rest}_{B / A} , \]
where $\opn{rest}_{B / A}$ is the forgetful functor. These are adjoint.
We get induced derived functors
\begin{equation} \label{equ:31}
 B \ot^{\mrm{L}}_A - : \
\til{\dcat{D}}(\cat{DGMod} A) \rightleftarrows \til{\dcat{D}}(\cat{DGMod} B) \
: \opn{rest}_{B / A} 
\end{equation}
that are also adjoint.
 
\begin{prop}
If $A \to B$ is a quasi-isomorphism, then the functors \tup{(\ref{equ:31})} are
equivalences.
\end{prop}

We say that $A$ is {\em strongly commutative} if 
$b \cd a = (-1)^{i + j} \cd a \cd b$ and $c^2 = 0$
for all $a \in A^i$, $b \in A^j$ and $c \in A^k$, where $k$ is odd. 
We call $A$ {\em nonpositive} if $A^i = 0$ for al $i > 0$. 

Let $f : A \to B$ be a homomorphism between nonpositive strongly commutative DG 
algebras. A K-flat DG algebra resolution of $B$ relative to $A$ is a 
factorization of 
$f$ into $A \xar{g} \til{B} \xar{h} B$, where $h$ is a quasi-isomorphism, and 
$\til{B}$ is a K-flat DG $A$-module. 
Such resolutions exist.
 
\begin{exa}
Take $A = \Z$ and $B := \Z / (6)$. We can take 
$\til{B}$ to be the Koszul complex
\[ \til{B} := ( \cdots 0 \to \Z \xar{6} \Z \to 0 \cdots ) \]
concentrated in degrees $-1$ and $0$. 
\end{exa}
 
\begin{exa}
For a homomorphism of commutative rings $A \to B$, the Hochschild cohomology of 
$B$ relative to $A$ is the cohomology of the complex 
\[ \opn{RHom}_{\til{B} \ot_A \til{B} }(B, B) , \]
where $\til{B}$ is a K-flat resolution as above. 
\end{exa}

\section{Commutative Dualizing Complexes}
 
I will talk about dualizing complexes over commutative rings.
There is a richer
theory for schemes, but there is not enough time for it. See \cite{Ha},
\cite{Ye2}, \cite{Ne1}, \cite{Ye5}, \cite{AJL}, \cite{LH} and their references.
 
Let $A$ be a noetherian commutative ring. 
We denote by $\dcat{D}_{\mrm{f}}^{\mrm{b}}(\cat{Mod} A)$
the subcategory of $\dcat{D}(\cat{Mod} A)$ consisting of bounded complexes whose
cohomologies are finitely generated $A$-modules. 
This is a full triangulated subcategory.
 
A complex $M \in \dcat{D}(\cat{Mod} A)$ is said to have {\em finite injective
dimension} if it has a bounded injective resolution. Namely there is a
quasi-isomorphism  $M \to I$ for some  bounded complex of injective $A$-modules
$I$. Note that such $I$ is a  K-injective complex. 

Take any $M \in \dcat{D} (\cat{Mod} A)$.
Because $A$ is commutative, we have a triangulated functor
\[ \opn{RHom}_A(-, M) : \dcat{D} (\cat{Mod} A)^{\mrm{op}} \to 
\dcat{D}(\cat{Mod} A) . \]
Cf.\ Example \ref{exa:5}.
 
\begin{dfn}
A {\em dualizing complex} over $A$ is a complex 
$R \in \dcat{D}_{\mrm{f}}^{\mrm{b}}(\cat{Mod} A)$ with finite injective
dimension, such that the canonical morphism 
\[ A \to \opn{RHom}_A(R, R) \]
in $\dcat{D}(\cat{Mod} A)$ is an isomorphism.
\end{dfn}
 
If we choose a bounded injective resolution $R \to I$, then there is an
isomorphism of triangulated functors 
\[ \opn{RHom}_A(-, R) \cong  \opn{Hom}_A(-, I) . \]

\begin{exa}
Assume $A$ is a {\em Gorenstein ring}, namely the free module $R := A$ has 
finite injective dimension.
There are plenty of Gorenstein rings; for instance any regular ring is
Gorenstein.
Then $R \in \dcat{D}_{\mrm{f}}^{\mrm{b}}(\cat{Mod} A)$,
and the reflexivity condition holds: 
\[ \opn{RHom}_A(R, R) \cong \opn{Hom}_A(A, A) \cong A . \]
We see that the module $R = A$ is a dualizing complex over the ring $A$. 
\end{exa}
 
Here are several important results from \cite{Ha}. 

\begin{thm}[Duality]
Suppose $R$ is a dualizing complex over $A$. Then the triangulated functor 
\[ \opn{RHom}_A(-, R) : 
\dcat{D}_{\mrm{f}}^{\mrm{b}} (\cat{Mod} A)^{\mrm{op}} \to
\dcat{D} _{\mrm{f}}^{\mrm{b}}(\cat{Mod} A)  \]
is an equivalence.
\end{thm}
 
\begin{thm}[Uniqueness]
Suppose $R$ and $R'$ are dualizing complexes over $A$, and $\opn{Spec} A$ is
connected. Then there is an invertible module $P$ and an integer $n$ such that
$R' \cong R \ot_A P[n]$
in $\dcat{D}_{\mrm{f}}^{\mrm{b}}(\cat{Mod} A)$. 
\end{thm}
 
\begin{thm}[Existence]
If $A$ has a dualizing complex, and $B$ is a finite type $A$-algebra, then 
$B$ has a dualizing complex. 
\end{thm}

\section{Noncommutative Dualizing Complexes}

In the last three sections of the paper we concentrate on noncommutative rings. 
Before going into the technicalities, here is a brief motivational preface. 

Recall that one of the important tools of commutative ring theory is 
{\em localization at prime ideals}. For instance, a noetherian local 
commutative ring $A$, with maximal ideal $\m$, is called a {\em regular local 
ring}  if
\[ \opn{dim} A = \opn{rank}_{A / \m}(\m / \m^2) . \]
(Here $\opn{dim}$ is Krull dimension.) A noetherian commutative ring $A$ is 
called {\em regular} if all its local rings $A_{\p}$ are regular local rings.

It is known that regularity can be described in homological terms. Indeed, if 
$\opn{dim} A < \infty$, then it is regular iff it has {\em finite global 
cohomological dimension}. Namely there is a natural number $d$, such that 
$\opn{Ext}^i_A(M, N) = 0$ for all $i > d$ and $M, N \in \cat{Mod} A$.  
 
Now consider a noetherian noncommutative ring $A$. (This is short for: $A$ is 
not-necessarily-commutative, and left-and-right noetherian.)
Localization of $A$ is almost never possible (for good reasons). 
A very useful substitute for localization (and other tools of commutative rings 
theory) is {\em noncommutative homological algebra}. By this we mean the study 
of the derived functors $\opn{RHom}_A(-, -)$, $\opn{RHom}_{A^{\mrm{op}}}(-, -)$
and $- \ot_A^{\mrm{L}} -$ of formulas (\ref{eqn:50}), 
(\ref{eqn:51}) and (\ref{eqn:52}) respectively. 
Here $A^{\mrm{op}}$ is the opposite ring (the same addition, but
multiplication is reversed).
The homological criterion of regularity from the commutative framework is  
made the definition of regularity in the noncommutative framework -- see 
Definition \ref{dfn:100} below. This definition is the point of departure of 
{\em noncommutative algebraic geometry} of M. Artin et.\ al.\ (see the survey 
paper \cite{SV}). 
A surprising amount of structure can be expressed in terms of 
noncommutative homological algebra. A few examples are sprinkled in the text, 
and many more are in the references.

\begin{dfn} \label{dfn:100}
A noncommutative ring $A$ is called {\em regular} if there is a natural number 
$d$, such that 
$\opn{Ext}^i_A(M, N) = 0$ and  $\opn{Ext}^i_{A^{\mrm{op}}}(M', N') = 0$
for all $i > d$, $M, N \in \cat{Mod} A$ and 
$M', N' \in \cat{Mod} A^{\mrm{op}}$.
\end{dfn}

For the rest of this section $A$ is a noncommutative noetherian ring. 
For technical reasons we assume that it is an algebra over a field $\K$. 
 
We denote by $A^{\mrm{e}} := A \ot_{\K} A^{\mrm{op}}$
the enveloping algebra.
Thus $\cat{Mod} A^{\mrm{op}}$ is the category of right $A$-modules,
and $\cat{Mod} A^{\mrm{e}}$ is the category of $\K$-central $A$-bimodules.
 
Any $M \in \cat{Mod} A^{\mrm{e}}$ gives rise to $\K$-linear 
functors 
\[ \opn{Hom}_A(-, M) : (\cat{Mod} A)^{\mrm{op}} \to 
\cat{Mod} A^{\mrm{op}}  \]
and 
\[ \opn{Hom}_{A^{\mrm{op}}}(-, M) : (\cat{Mod} A^{\mrm{op}})^{\mrm{op}} \to 
\cat{Mod} A . \]
These functors can be right derived, yielding  $\K$-linear triangulated 
functors 
\begin{equation} \label{eqn:50}
\opn{RHom}_A(-, M) : \dcat{D}(\cat{Mod} A)^{\mrm{op}} \to 
\dcat{D}(\cat{Mod} A^{\mrm{op}})
\end{equation}
and 
\begin{equation} \label{eqn:51}
\opn{RHom}_{A^{\mrm{op}}}(-, M) : 
\dcat{D}(\cat{Mod} A^{\mrm{op}})^{\mrm{op}}
\to \dcat{D}(\cat{Mod} A) .
\end{equation}
One way to construct these derived functors is to choose a 
K-injective resolution
$M \to I$ in $\dcat{K}(\cat{Mod} A^{\mrm{e}})$. 
Then (because $A$ is flat over $\K$) the complex $I$ is 
K-injective over $A$ and over $A^{\mrm{op}}$, and we get
\[ \opn{RHom}_A(-, M) \cong \opn{Hom}_A(-, I) \]
and
\[ \opn{RHom}_{A^{\mrm{op}}} (-, M)  \cong 
\opn{Hom}_{A^{\mrm{op}}}(-, I) . \] 

Note that even if $A$ is commutative, this setup is still meaningful -- not
all $A$-bimodules are $A$-central!
 
\begin{dfn}[\cite{Ye1}]
 A {\em noncommutative dualizing complex} over $A$ is a complex 
$R \in \dcat{D}^{\mrm{b}}(\cat{Mod} A^{\mrm{e}})$ satisfying these three
conditions:  
\begin{enumerate} 
\rmitem{i} The cohomology modules $\mrm{H}^q(R)$ are finitely generated 
over $A$ and over $A^{\mrm{op}}$.
 
\rmitem{ii} The complex $R$ has finite injective dimension over $A$ and over 
$A^{\mrm{op}}$.

\rmitem{iii} The canonical morphisms 
\[ A \to \opn{RHom}_A(R, R) \]
and 
\[ A \to \opn{RHom}_{A^{\mrm{op}}}(R, R) \]
in $\dcat{D}(\cat{Mod} A^{\mrm{e}})$ are isomorphisms.
\end{enumerate}
\end{dfn}

Condition (ii) implies that $R$ has a ``bounded bi-injective resolution'',
namely there is a quasi-isomorphism $R \to I$ 
in $\dcat{K}(\cat{Mod} A^{\mrm{e}})$,
with $I$ a bounded complex of bimodules that are injective on both sides.
 
\begin{thm}[Duality, \cite{Ye1}]
Suppose $R$ is a noncommutative dualizing complex over $A$. Then the
triangulated functor 
\[ \opn{RHom}_A(-, R) : 
\dcat{D}_{\mrm{f}}^{\mrm{b}} (\cat{Mod} A)^{\mrm{op}} \to
\dcat{D}_{\mrm{f}}^{\mrm{b}} (\cat{Mod} A^{\mrm{op}}) \]
is an equivalence, with quasi-inverse 
$\opn{RHom}_{A^{\mrm{op}}}(-, R)$.
\end{thm}
 
Existence and uniqueness are much more complicated than in the noncommutative
case. I will talk about them later.

\begin{exa}
The noncommutative ring $A$ is called {\em Gorenstein} if the bimodule 
$A$ has finite injective dimension on both sides. 
It is not hard to see that $A$ is Gorenstein iff it has a noncommutative
dualizing complex of the form $P[n]$, for some integer $n$ and {\em
invertible bimodule} $P$.
Here invertible bimodule is in the sense of Morita theory, namely there is
another bimodule $P^{\vee}$ such that 
\[ P \ot_A P^{\vee} \cong P^{\vee} \ot_A P \cong A \]
in $\cat{Mod} A^{\mrm{e}}$.
Any regular ring is Gorenstein.
\end{exa}
 
For more results about noncommutative Gorenstein rings see \cite{Jo} and
\cite{JZ}.

\section{Tilting Complexes and Derived Morita Theory}
 
Let $A$ and $B$ be noncommutative algebras over a field $\K$.  
Suppose $M \in \lb \dcat{D}(\cat{Mod} A \ot_{\K} B^{\mrm{op}})$
and $N \in \dcat{D}(\cat{Mod} B \ot_{\K} A^{\mrm{op}})$.
The left derived tensor product 
\begin{equation} \label{eqn:52}
M \ot_B^{\mrm{L}} N \in \dcat{D}(\cat{Mod} A \ot_{\K} A^{\mrm{op}})
\end{equation}
exists. It can be constructed by choosing a resolution $P \to M$ in 
$\dcat{K}(\cat{Mod} A \ot_{\K} B^{\mrm{op}})$, where $P$ is
a complex that's K-projective over $B^{\mrm{op}}$;
or by choosing a resolution $Q \to N$ in 
$\dcat{K}(\cat{Mod} B \ot_{\K} A^{\mrm{op}})$, where $Q$ is
a complex that's K-projective over $B$.

Here is a definition generalizing the notion of invertible bimodule. It is due
to Rickard \cite{Ri1}, \cite{Ri2}. 
 
\begin{dfn}
A complex 
\[ T \in \dcat{D}(\cat{Mod} A \ot_{\K} B^{\mrm{op}}) \]
is called a {\em two-sided tilting complex} over $A$-$B$ 
 if there exists a complex 
\[ T^{\vee} \in \dcat{D}(\cat{Mod} B \ot_{\K} A^{\mrm{op}}) \]
 such that 
\[ T \ot_B^{\mrm{L}} T^{\vee} \cong A \]
in $\dcat{D}(\cat{Mod} A^{\mrm{e}})$, and 
\[ T^{\vee} \ot_A^{\mrm{L}} T \cong B \]
in $\dcat{D}(\cat{Mod} B^{\mrm{e}})$.
 
When $B = A$ we say that $T$ is a two-sided tilting complex over $A$.
\end{dfn}

The complex $T^{\vee}$ is called a quasi-inverse of $T$. It is unique up to
isomorphism in $\dcat{D}(\cat{Mod} B \ot_{\K} A^{\mrm{op}})$. Indeed we have
this result:
 
\begin{prop}
Let $T$ be a two-sided tilting complex.
 
\begin{enumerate} 
\item The quasi-inverse $T^{\vee}$ is isomorphic to 
$\opn{RHom}_A(T, A)$.
 
\item $T$ has a bounded bi-projective resolution $P \to T$.
\end{enumerate}
\end{prop}
 
\begin{dfn}
The algebras $A$ and $B$ are said to be {\em derived Morita equivalent} if
there is a $\K$-linear triangulated equivalence 
\[ \dcat{D}(\cat{Mod} A) \approx \dcat{D}(\cat{Mod} B) . \]
\end{dfn}
 
\begin{thm}[\cite{Ri2}]
The $\K$-algebras $A$ and $B$ are derived Morita equivalent iff there exists a
two-sided tilting complex over $A$-$B$.
\end{thm} 

Here is a result relating dualizing complexes and tilting complexes. 
 
\begin{thm}[Uniqueness, \cite{Ye3}] \label{thm:7}
Suppose $R$ and $R'$ are noncommutative dualizing complexes over $A$.
Then the complex 
\[ T := \opn{RHom}_A(R, R') \]
is a two-sided tilting complex over $A$, and 
\[ R' \cong R \ot_A^{\mrm{L}} T \]
in $\dcat{D}(\cat{Mod} A^{\mrm{e}})$. 
\end{thm}

It is easy to see that if $T_1$ and $T_2$ are two-sided tilting complexes over
$A$, then so is $T_1 \ot_A^{\mrm{L}} T_2$. 
This leads to the next definition. 
 
\begin{dfn}[\cite{Ye3}]
Let $A$ be a noncommutative $\K$-algebra.
The {\em derived Picard group} of $A$ is the group $\opn{DPic}_{\K}(A)$
whose elements are the isomorphism classes (in 
$\dcat{D}(\cat{Mod} A^{\mrm{e}})$)
of two-sided tilting complexes.
The multiplication is induced by the operation
$T_1 \ot_A^{\mrm{L}} T_2$, 
and the identity element is the class of $A$.
\end{dfn}

Here is a consequence of Theorem \ref{thm:7}.
 
\begin{cor}
Suppose the noncommutative $\K$-algebra $A$ has at least one dualizing complex. 
Then operation $R \ot_A^{\mrm{L}} T$
induces a simply transitive right action of the group 
$\opn{DPic}_{\K}(A)$ on the set of isomorphism classes of dualizing
complexes.
\end{cor}
 
It is natural to ask about the structure of the group $\opn{DPic}(A)$.

\begin{thm}[\cite{RZ}, \cite{Ye3}]
If the ring $A$ is either commutative \tup{(}with nonempty connected 
spectrum\tup{)} or local, then 
\[ \opn{DPic}_{\K}(A) \cong \opn{Pic}_{\K}(A) \times \Z . \]
\end{thm}
 
Here $\opn{Pic}_{\K}(A)$ is the noncommutative Picard group of $A$, made up of
invertible bimodules. If $A$ is commutative, then 
\[ \opn{Pic}_{\K}(A) \cong \opn{Aut}_{\K}(A) \ltimes \opn{Pic}_{A}(A) , \]
where $\opn{Pic}_{A}(A)$ is the usual (commutative) Picard group of $A$. 
A noncommutative ring $A$ is said to be local if $A / \mfrak{r}$ is a simple 
artinian ring, where $\mfrak{r}$ is the Jacobson radical.
 
For nonlocal noncommutative rings the group $\opn{DPic}_{\K}(A)$ is bigger. 
See the paper \cite{MY} for some calculations. 
These calculations are related to CY-dimensions of some rings; cf.\ Example 
\ref{exa:10}.

\section{Rigid Dualizing Complexes}
 
The material in this final section is largely due to Van den Bergh
\cite{VdB1}. His results were extended by J. Zhang and myself. 
Again $A$ is a noetherian noncommutative algebra over a field $\K$, and 
$A^{\mrm{e}} = A \ot_{\K} A^{\mrm{op}}$. 
 
Take $M \in \cat{Mod} A^{\mrm{e}}$. Then the $\K$-module
$M \ot_{\K} M$ has four commuting actions by $A$, which we arrange
as follows.
The algebra $A^{\mrm{e};\, \mrm{in}} := A^{\mrm{e}}$ acts
on $M \ot_{\K} M$ by 
\[ (a_1 \ot a_2) \cdot_{\tup{in}} (m_1 \ot m_2) :=
 (m_1 \cdot a_2) \ot (a_1 \cdot m_2) , \]
and the algebra $A^{\mrm{e};\, \mrm{out}} := A^{\mrm{e}}$ acts by  
\[ (a_1 \ot a_2) \cdot_{\tup{out}}  (m_1 \ot m_2) := 
(a_1 \cdot m_1) \ot (m_2 \cdot a_2) . \]
The bimodule $A$ is viewed as an object of $\dcat{D}(\cat{Mod} A^{\mrm{e}})$
in the obvious way. 
 
Now take $M \in \dcat{D}(\cat{Mod} A^{\mrm{e}})$. We define the {\em square}
of $M$ to be the complex 
\[ \opn{Sq}_{A / \K}(M) := 
\opn{RHom}_{A^{\mrm{e}; \, \mrm{out}}}(A, M \ot_{\K} M) 
\in \dcat{D}(\cat{Mod} A^{\mrm{e}; \, \mrm{in}}) . \]
We get a functor 
\[ \opn{Sq}_{A / \K} : \dcat{D}(\cat{Mod} A^{\mrm{e}}) \to 
\dcat{D}(\cat{Mod} A^{\mrm{e}})  . \]

This is not an additive functor. Indeed, it is a quadratic functor: given an
element $a \in \mrm{Z}(A)$ and a morphism 
$\phi : M \to N$ in $\dcat{D}(\cat{Mod} A^{\mrm{e}})$,
one has 
\[ \opn{Sq}_{A / \K}(a \cd \phi) = \opn{Sq}_{A / \K}(\phi \cd a) = 
a^2 \cd \opn{Sq}_{A / \K}(\phi) . \]

Note that the cohomologies of $\opn{Sq}_{A / \K}(M)$ are 
\[ \mrm{H}^j(\opn{Sq}_{A / \K}(M)) = 
\opn{Ext}^j_{A^{\mrm{e}}}(A, M \ot_{\K} M) , \]
so they are precisely the Hochschild cohomologies of $M \ot_{\K} M$.
 
A {\em rigid complex} over $A$ (relative to $\K$) is a pair $(M, \rho)$ 
consisting of a complex $M \in \dcat{D}(\cat{Mod} A^{\mrm{e}})$,
and an isomorphism 
\[ \rho : M \iso \opn{Sq}_{A / \K}(M) \]
in $\dcat{D}(\cat{Mod} A^{\mrm{e}})$.
 
Let  $(M, \rho)$  and $(N, \sigma)$ be rigid complexes over $A$.
A {\em rigid morphism}
\[ \phi :  (M, \rho) \to (N, \sigma) \]
is a morphism $\phi : M \to N$ in $\dcat{D}(\cat{Mod} A^{\mrm{e}})$, such that
the diagram
\[ \UseTips  \xymatrix @C=6ex @R=5ex {
M 
\ar[r]^(0.4){\rho}
\ar[d]_{\phi}
&
\opn{Sq}_{A / \K}(M)
\ar[d]^{\opn{Sq}_{A / \K}(\phi)}
\\
N
\ar[r]^(0.4){\sigma}
&
\opn{Sq}_{A / \K}(N)
} \]
is commutative. 

\begin{dfn}[\cite{VdB1}]
A {\em rigid dualizing complex} over $A$ (relative to $\K$) is a rigid
complex $(R, \rho)$ such that $R$ is a dualizing complex. 
\end{dfn}

\begin{thm}[Uniqueness, \cite{VdB1}, \cite{Ye3}]
Suppose $(R, \rho)$ and $(R', \rho')$ are both rigid dualizing complexes over
$A$. Then there is a unique rigid isomorphism 
\[  \phi :  (R, \rho) \iso (R', \rho') . \]
\end{thm}

As for existence, let me first give an easy case. 
 
\begin{prop}
If $A$ is finite over its center $\opn{Z}(A)$, and $\opn{Z}(A)$
is finitely generated as 
$\K$-algebra, then $A$ has a rigid dualizing complex. 
\end{prop}
 
Actually, in this case it is quite easy to write down a formula for the rigid
dualizing complex.
 
In the next existence result, by a filtration 
$F = \{ F_i(A) \}_{i \in \N}$ of the
algebra $A$ we mean an ascending exhaustive filtration
by $\K$-submodules, such that $1 \in F_0(A)$ and 
$F_i(A) \cd F_j(A) \subset F_{i + j}(A)$.
Such a filtration gives rise to a graded $\K$-algebra 
\[ \opn{gr}^F(A) = \bigoplus_{i \geq 0}\ \opn{gr}^F_i(A) . \]
 
\begin{thm}[Existence, \cite{VdB1}, \cite{YZ3}]
Suppose $A$ admits a filtration $F$, such that $\opn{gr}^F(A)$ is 
finite over its center $\opn{Z}(\opn{gr}^F(A))$, 
and $\opn{Z}(\opn{gr}^F(A))$ is finitely generated as $\K$-algebra. Then $A$ has
a rigid dualizing complex. 
\end{thm}
 
This theorem applies to the ring of differential operators 
$\mcal{D}(C)$, where $C$ is a smooth commutative $\K$-algebra 
(and $\opn{char} \K = 0$). 
It also applies to any quotient of 
the universal enveloping algebra $\opn{U}(\g)$ of a finite dimensional Lie
algebra $\g$. 
 
I will finish with some examples.

\begin{exa}
Let $A$ be a noetherian $\K$-algebra satisfying these two conditions:
\begin{itemize}
\item $A$ is smooth, namely the $A^{\mrm{e}}$-module $A$ has finite projective
dimension.
\item There is an integer $n$ such that 
\[ \opn{Ext}^j_{A^{\mrm{e}}}(A, A^{\mrm{e}}) \cong 
\begin{cases}
A & \tup{if} \ j = n
\\
0 & \tup{otherwise}.
\end{cases}
\]
\end{itemize}
Then $A$ is a regular ring (Definition \ref{dfn:100}), and the complex 
$R := A[n]$
is a rigid dualizing complex over $A$. 
Such an algebra $A$ is called an {\em $n$-dimensional Artin-Schelter regular
algebra}, or an {\em $n$-dimensional Calabi-Yau algebra}. 
\end{exa}

\begin{exa}
Let $\g$ be an $n$-dimensional Lie algebra, and $A := \opn{U}(\g)$, the
universal enveloping algebra.
Then the rigid dualizing complex of $A$ is 
$R := A^{\sigma}[n]$,
where $A^{\sigma}$ is the trivial bimodule $A$, 
twisted on the right by an automorphism $\sigma$.
Using the Hopf structure of $A$ we can express $A^{\sigma}$ like this:
\[ A^{\sigma} \cong \opn{U}(\g) \ot_{\K} \bwedge^n (\g) , \]
the twist by the $1$-dimensional representation $\bwedge^n (\g)$.
See \cite{Ye4}. 
So $A$ is a {\em twisted Calabi-Yau} algebra.
If $\g$ is semi-simple then there is no twist, and $A$ is Calabi-Yau. This was
used by Van den Bergh in his duality for Hochschild (co)homology \cite{VdB2}.
\end{exa}

\begin{exa} \label{exa:10}
Let 
\[ A := \bmat{\K & \K \\ 0 & \K} \]
the $2 \times 2$ matrix algebra. 
The rigid dualizing complex here is 
\[ R := \opn{Hom_{\K}}(A, \K) . \]
It is known that 
\[ R \ot_A^{\mrm{L}} R \ot_A^{\mrm{L}} R \cong A[1] \]
in $\dcat{D}(\cat{Mod} A^{\mrm{e}})$.
So $A$ is a {\em Calabi-Yau algebra of dimension $\frac{1}{3}$}.
See \cite{Ye3}, \cite{MY}. 
\end{exa}

%\section{References}

\end{document}